\documentclass[12pt]{amsart}
\usepackage{amscd, amssymb,latexsym,amsmath, amscd, amsmath}
\usepackage{tikz-cd}
\usepackage{extarrows}
\usepackage{mathtools}
\usepackage{hyperref}
\usepackage[backend=biber,
style=alphabetic]{biblatex}
\hypersetup{%
        colorlinks,
        breaklinks=true,
        plainpages=false,%
        citecolor=blue,
        linkcolor=blue,
        urlcolor=blue,
        bookmarksopen=true,%
        bookmarksnumbered=false,%
        bookmarksdepth=5%
}

\newcommand{\Q}{\mathbb{Q}}

\newcommand{\Z}{\mathbb{Z}}

\newcommand{\F}{\mathbb{F}}

\newcommand{\Ext}{{\rm Ext}}

\newcommand{\Fp}{\mathbb{F}_p}
\newcommand{\fr}{\mathrm{Fr}}
\newcommand{\Fq}{\mathbb{F}_q}
\newcommand{\im}{\rm{Im}\>}
\newcommand{\Ker}{\rm{Ker}\>}

\newcommand\restr[2]{{% we make the whole thing an ordinary symbol
  \left.\kern-\nulldelimiterspace % automatically resize the bar with \right
  #1 % the function
  \littletaller % pretend it's a little taller at normal size
  \right|_{#2} % this is the delimiter
  }}
\newcommand{\littletaller}{\mathchoice{\vphantom{\big|}}{}{}{}}

\def\GL{{\rm GL}}

\newcommand{\IZind}{\mathrm{ind}_{IZ}^G\>}
\newcommand{\KZind}{\mathrm{ind}_{KZ}^G\>}

\newcommand{\id}{\mathrm{id}}
\newcommand{\Id}{\mathrm{Id}}

\begin{document}

\newtheorem{theorem}{Theorem}
\newtheorem{thm}[equation]{Theorem}
\newtheorem{prop}{Proposition}
\newtheorem{cor}[equation]{Corollary}
\newtheorem{conj}{Conjecture}
\newtheorem{lemma}{Lemma}
\newtheorem{corollary}{Corollary}
\newtheorem{question}[equation]{Question}

\newtheorem{conjecture}[theorem]{Conjecture}

\newtheorem{example}[equation]{Example}
\numberwithin{equation}{section}

\newtheorem{remark}{Remark}

\title{Iwahori-Hecke model for the universal supersingular representation}

\author[Anand Chitrao]{Anand Chitrao}
\author[Asfak Soneji]{Asfak Soneji}
\address{Department of Mathematics, HRI, Prayagraj (Allahabad) - 211019, India}
\email{anandchitrao@gmail.com}
\address{Department of Basic Sciences, IITRAM, Ahmedabad - 380026, India}
\email{sonejiasfak96@gmail.com}

\begin{abstract}
Let $F$ be a non-archimedean local field with residue field $\F_q$. When $F$ is a finite extension of $\mathbb{Q}_{p}$, Anandavardhanan-Borisagar and Anandavardhanan-Jana introduced an Iwahori-Hecke model for the universal supersingular representation $\pi_r$ in the \emph{regular} case $ 0 < r < q-1$. When $F=\mathbb{Q}_{p}$, the first author introduced an Iwahori-Hecke model for $\pi_r$ when $r = 0, \ p - 1$.
We extend this result to an arbitrary local field $F$ for $r = 0, \ q - 1$. 

We also write down an explicit non-split self-extension of $\pi_r$ which has a four-dimensional space of $I(1)$-invariants when $F = \mathbb{Q}_p$ and $r = 0, \ p-1$.
\end{abstract}
%%% ----------------------------------------------------------------------
\maketitle
%%% ----------------------------------------------------------------------

%---------------------------------------------------------------------------

\section{Notation and Conventions}
For a prime $p$, let $\F_q=\F_{p^{f}}$ be a degree $f$ extension of the finite field $\F_p$. Let $\mathbb{\overline{F}}_{p}$ denote an algebraic closure of $\F_p$. We fix an embedding of $\F_q \hookrightarrow \mathbb{\overline{F}}_{p}$. Let $F$ be a local field with ring of integers $\mathcal{O}$, uniformizer $\varpi$ and residue field $\Fq$. Let $G$  be the group ${\rm{GL}_2}(F)$ and $Z$ be its center. Let $I_1=\left\{[{\mu}]\mid
{\mu}\in \mathbb{F}_q\right\} \subset \mathcal{O}$,
where for $\mu \in \F_q$, we denote its Teichm\"uller lift in $\mathcal O$ by $[\mu]$. Let $K= \rm{GL_2}(\mathcal{O})$ be the maximal compact subgroup. Let $I$ denote the Iwahori subgroup of upper triangular matrices mod $\varpi$, and $I(1)$ be the subgroup of $K$ consisting of matrices that are upper triangular unipotent mod $\varpi$. Let 
$$w :=\begin{pmatrix}
    0 & 1 \\
    1 & 0
\end{pmatrix}, \alpha := \begin{pmatrix}
    1 & 0 \\
    0 & \varpi
\end{pmatrix} \  \text{and} \ \ \beta :=\begin{pmatrix}
    0 & 1 \\
    \varpi & 0
\end{pmatrix} = \alpha w.$$
For an integer $r\geq 0$, we write $r=r_0+r_1p +\cdots+ r_{f-1}p^{f-1}$ where $r_j \geq 0$ for $0 \leq j \leq f-1$. Let $V_{r_j}$ denote the symmetric $r_j$-th power of the standard two-dimensional representation of $\rm{GL}_{2}(\F_q)$. We identify it with the homogeneous polynomials of degree $r_{j}$ in variables $X_j$ and $Y_j$ for all $0 \leq j \leq f-1$. By $V_{r_j}^{\fr^{j}}$ we mean that the action of $\rm{GL}_{2}(\F_q)$ is twisted by the $j$-th power of Frobenius. We view these as representations of $K$, via the reduction mod $\varpi$ map and then as representations of $KZ$ by making the scalar matrix $\varpi$ act trivially. For $r = r_0 + r_1 p + \cdots + r_{f - 1}p^{f - 1}$, let $V_r$ be the representation $V_{r_0} \otimes_{\Fq} V_{r_1}^{\fr} \otimes_{\Fq} \cdots \otimes_{\Fq} V_{r_{f - 1}}^{\fr^{f - 1}}$ of $KZ$. For $0 \leq r < q - 1$, let $d^r$ be the character $ I \to \Fq^{\times}$ sending $ \begin{pmatrix}a & b \\ c & d\end{pmatrix}$ to $d^r \mod \varpi$ and extend it to $IZ$ by mapping the scalar matrix $\varpi \in Z $ to $1$. If $0 < r < q - 1$, we say that $d^r$ is \emph{regular} and if $r = 0$, we say that $d^r$ is \emph{non-regular}.

Let $\KZind \sigma$ denote the compact induction of a representation $\sigma$ of $KZ$. Following notation of Barthel-Livne \cite{BL94}, for $g \in G$ and $v \in \sigma$, the symbol $[g, v] \in \KZind \sigma$ will be used to denote the function from $G$ to $\sigma$ which is supported on the coset $KZg^{-1}$ and which sends $g^{-1}$ to $v$. To distinguish elements of different induced representations, we will denote elements of compact induction  $\operatorname{ind}_{IZ}^{G}\chi$ by $[[g, v]]$ and those of $\operatorname{ind}_{IZ}^{KZ} \sigma$ by $[[[g,v]]]$.
For an irreducible representation $\chi$, we define the Iwahori-Hecke algebra as the endomorphism algebra $\mathrm{End}_G(\IZind \chi)$. Similarly, for an irreducible representation $\sigma$ of $KZ$, the spherical Hecke algebra is defined as $\mathrm{End}_G(\KZind \sigma)$. For $j \in \Z/f\Z$ let $\theta_{j}=X_{j}Y_{j-1}^{p}-Y_{j}X_{j-1}^{p}$ be the Ghate-Jana polynomials defined in \cite{GJ25}.

\section{Introduction}
For a local field $F$, the classification of irreducible smooth mod $p$ representations admitting a central character of $G$ began with the notable work of Barthel and Livn\'e \cite{BL94}. In addition to many interesting results, they proved the existence of irreducible smooth representations, known as supersingular representations, which cannot be obtained as a subquotient of a parabolically induced representation.

It is proved in \cite{BL94} that up to a twist by a character, a supersingular representation can be realized as an irreducible quotient of a universal representation $\pi_{r}=\frac{\KZind V_r}{(T)}$ for $0 \leq r \leq q - 1$, where the spherical Hecke operator $T$ is such that ${\rm End}_G\left(\operatorname{ind}_{KZ}^{G} V_r\right) \cong \Fq[T]$. Computing $I(1)$-invariants, Breuil showed that $\pi_r$ is irreducible \cite{Bre03} when $F=\mathbb{Q}_{p}$. However, if $[F:\Q_p] > 1$, the structure of $\pi_r$ is mysterious. For instance, they are of infinite length. In \cite{Sch11} and \cite{Sch14}, Schein studied the $K$-structure of the $G$-representations $\pi_r$ when $[F: \Q_p] > 1$ and established a criterion for irreducibility concerning supersingular representations. It can be seen in \cite{BP12} that explicit construction of a supersingular representation is a difficult task when $F \supsetneq \Q_p$ is unramified. Schein \cite{Sch23} constructed a family of supersingular representations with good socle for a ramified extension $F$ of $ \mathbb{Q}_{p}$, with the ramification index $e \leq \frac{p-1}{2}$ using diagrams. For a finite extension $F$ of $\mathbb{Q}_{p}$, Hendel \cite{Hen19} computed the space of $I(1)$-invariants of $\pi_r$, and when $F$ is unramified over $\Q_p$, constructed a universal supersingular representation having all supersingular representations with a good $K$-socle as quotients. Jana has extended Hendel's work concerning the $I(1)$-invariants of $\pi_r$ in \cite{Jan23}.

A universal supersingular representation can also be constructed using the Iwahori-Hecke operators instead of the spherical Hecke operator $T$  for a finite extension of $\mathbb{Q}_{p}$ (see \cite{AB15}, \cite{AJ21}). For $0<r<q-1$, consider the \emph{regular} character $d^r$ of $IZ$. It is shown in \cite{BL94} that the corresponding Iwahori-Hecke algebra is commutative and is given by
$${\rm End}_G\left(\operatorname{ind}_{IZ}^{G} d^r\right) \cong \frac{\Fq\left[T_{-1,0},T_{1,2}\right]}{(T_{-1,0}T_{1,2}, T_{1,2}T_{-1,0})}.$$ Here, $T_{-1,0}$ and $T_{1,2}$ are Iwahori-Hecke operators as defined in \cite{AJ21}.
When $F$ is a totally ramified extension of $\mathbb{Q}_{p}$, \cite[Remark 1]{AB15}, we have
\begin{eqnarray}\label{totally ramified case}
   {\im} T_{-1,0}= {\Ker} T_{1,2} \quad \text{and} \quad {\im} T_{1,2}= {\Ker} T_{-1,0}.
\end{eqnarray}

When $F$ is not totally ramified, there are strict inclusions 
\begin{eqnarray}\label{not totally ramified case}
   {\im} T_{-1,0} \subsetneq {\Ker} T_{1,2} \quad \text{and} \quad {\im} T_{1,2} \subsetneq {\Ker} T_{-1,0} 
\end{eqnarray}
as shown in \cite[Remark $5$]{AJ21}. When $F$ is a finite extension of $\mathbb{Q}_p$, the comparison between Iwahori induction and spherical induction is given in \cite[Theorem 1.1]{AJ21}.

The Iwahori model for the universal supersingular representations is also useful when one wants to compute the invariants in $\pi_r$ under the action of principal congruence subgroups. This is because using \eqref{totally ramified case} one can write $\pi_r$ as a quotient of an induced representation by kernels of certain Iwahori-Hecke operators. This idea is worked out in \cite[Theorem 1.2]{AB13} to give an easier proof of Morra's result \cite[Theorem 1.4]{Mor13} on the $K$-structure of the invariants of $\pi_r$ under the principal congruence subgroups when $r \neq 0, p - 1$.

For the \emph{non-regular} case $r\in \{0, q-1\}$, the Iwahori-Hecke algebra is non-commutative and is generated by Iwahori-Hecke operators $T_{1,0}$ and $T_{1,2}$ satisfying 
\begin{eqnarray}\label{relations in non-reqular case}
    T_{1,0}^2=\Id \quad \text{and} \quad T_{1,2}T_{1,0}T_{1,2}=-T_{1,2}
\end{eqnarray}
as given in \cite[Proposition 11]{BL94}.
As far as the remaining cases are concerned, the first author \cite[Theorem 3.13]{Chi25} gives the Iwahori-Hecke model for the \emph{non-regular} case but exclusively for $\mathbb{Q}_p$, i.e., for $ r=0$ and $p-1$. In this article, we extend this model to include arbitrary local fields $F$, see Theorem \ref{main theorem}. We apply the Iwahori-Hecke model for ${\rm{GL}_2}(\mathbb{Q}_{p})$ to construct self-extensions of supersingular representations, which have four-dimensional $I(1)$-invariants in the \emph{non-regular} case, as is outlined in Theorems \ref{T3} and \ref{T4}. We prove Lemmas \ref{L2} and \ref{L1}, 
which, when compared with \eqref{not totally ramified case}, show the striking dissimilarity between the \emph{regular} and the \emph{non-regular} cases.

One of the motivations behind this paper was to write down an Iwahori theoretic modular local Langlands correspondence between $\mathrm{Gal}(\overline{F}/F)$ and $\GL_2(F)$ as was done in \cite{Chi25} for $F = \Q_p$. The correspondence in the language of induction from the Iwahori subgroup (times the center) is especially useful while computing the reduction of semi-stable representations of $\mathrm{Gal}(\overline{F}/F)$ when $F = \Q_p$. It turns out that the automorphic representation associated by the $p$-adic local Langlands correspondence with a two-dimensional, irreducible, semi-stable representation of the absolute Galois group of $\Q_p$ has a dense sub-representation isomorphic to a quotient of a representation induced from the Iwahori subgroup to $\GL_2(\Q_p)$. Therefore, when computing the reduction mod $p$ of the automorphic representation, one must work with Iwahori induction. Using this idea, it is possible to compute the reduction mod $p$ of two-dimensional, irreducible, semi-stable representations of the absolute Galois group of $\Q_p$ with \emph{arbitrarily large} Hodge-Tate weights. This has been done in \cite{CG23} for weights at most $p + 1$, including the previously untreated weights $p$ and $p + 1$, and is further explored in \cite{CG25} for weights at most $p^2$.

\section{preliminaries}
We state Lucas' classical result in modular combinatorics, which computes the residue of the binomial coefficient ${n \choose i}$  
modulo $p$.
\begin{theorem}[Lucas]\label{lucas}
Let $n,i\in \mathbb{N}$ be such that $n=\sum\limits_{\substack{j=0}}^kn_jp^j$ and $i=\sum\limits_{\substack{j=0}}^ki_jp^j$, where $0\leq n_j\leq p-1$ and $0\leq i_j\leq p-1$. Then \[{n\choose i}\equiv\prod_{j=0}^k{n_j\choose i_j}\mod p.\] 
\end{theorem}

\begin{corollary}\label{prime divides bio coeff}
Let $n, i\in \mathbb{N}$. Then $p$ divides ${n\choose i}$ if and only if $n_j<i_j$ for some $0\leq j\leq k$.
\end{corollary}
As mentioned in the introduction, the Iwahori-Hecke algebra is non-commutative in the \emph{non-regular} case. It is generated by the Iwahori-Hecke operators $T_{1,0}$ and $T_{1,2}$, which satisfy the relations as stated in \eqref{relations in non-reqular case}. We write down the explicit formulas for these operators as given in \cite{BL95}:  
\begin{eqnarray}\label{Formulas T10 and T12 }
        T_{1, 0}[[g, 1]] = [[g\beta, 1]] \text{ and } T_{1, 2}[[g, 1]] = \sum_{\lambda \in I_1}\left[\left[g\begin{pmatrix}1 & 0 \\ \varpi\lambda & \varpi\end{pmatrix}, 1\right]\right].
    \end{eqnarray}
 According to \cite[Lemma 8(4)]{BL95}, the operator $T_{-1,0}$ satisfies the following relation:   \begin{eqnarray}\label{T-10 in the non-commutative case}
        T_{-1, 0} = T_{1, 0}T_{1, 2}T_{1, 0}.
    \end{eqnarray}
Its formula is given as follows:
\begin{eqnarray}\label{Formulas for T-10 and T12}
T_{-1, 0}[[g, v]] = \sum_{\lambda \in I_1}\left[\left[g\begin{pmatrix}\varpi & \lambda \\ 0 & 1\end{pmatrix}, v\right]\right].
\end{eqnarray}
The formula for the spherical Hecke operator $T$ is given by
\begin{eqnarray}\label{Formula for T in the case r=0}
\hspace{1cm} T[g, 1] = \sum_{\lambda \in I_1}\left[g\begin{pmatrix}\varpi & \lambda \\ 0 & 1\end{pmatrix}, 1\right] + \left[g\begin{pmatrix}1 & 0 \\ 0 & \varpi\end{pmatrix}, 1\right], & \textit{ if }  r=0.  \end{eqnarray}
\begin{eqnarray}\label{Formula for T in the case r=q-1}
\hspace{1cm} T[g, X_0^{p-1} \cdots X_{f-1}^{p-1}] = \sum_{\lambda \in I_1}\left[g\begin{pmatrix}\varpi & \lambda \\ 0 & 1\end{pmatrix}, X_0^{p-1} \cdots X_{f-1}^{p-1}\right], & \textit{ if }  r=q-1.  \end{eqnarray}

\section{Comparison between Iwahori induction and spherical induction}

In this section, we prove the main comparison theorem between Iwahori induction and spherical induction.

The following lemma extends a case of \cite[Lemma 2.15]{GV22}. This lemma also appears in the ongoing work \cite{GP25}.

\begin{lemma}\label{Sum of products of binomial coeffs is a sign}
For $0 \leq j \leq f - 1$, fix integers $r_j \geq 0$ such that $0<r = r_0 +  r_1p + \ldots + r_{f-1}p^{f - 1}$ is divisible by $q - 1$. Then for any integer $i = i_0 +  i_1p + \ldots + i_{f - 1}p^{f - 1} \in [0, q - 1)$, we have
    \[
        \sum_{\substack{(k_0, \ldots, k_{f - 1}) = (0, \ldots, 0) \\ k \equiv i \!\! \mod q - 1}}^{(r_0, \ldots, r_{f - 1})} {r_0 \choose k_0}\cdots {r_{f - 1} \choose k_{f - 1}} \equiv 
        \begin{cases}
            (-1)^i \mod p & \text{ if } i \in (0, q - 1), \\
            2 \mod p& \text{ if } i = 0.
        \end{cases}
    \]
\end{lemma}
\begin{proof}
    Consider the polynomial $f(x) = (1 + x)^{r_0}(1 + x^p)^{r_1} \cdots (1 + x^{p^{f - 1}})^{r_{f - 1}}$ in the ring $R = \Fp[x]/(x^{q - 1} - 1)$. For $i \in [0, q - 1)$, the coefficient of $x^i$ in $f(x)$ is
    \[
        \sum_{\substack{(k_0, \ldots, k_{f - 1}) = (0, \ldots, 0) \\ k \equiv i \!\! \mod q - 1}}^{(r_0, \ldots, r_{f - 1})} {r_0 \choose k_0}\cdots {r_{f - 1} \choose k_{f - 1}}.
    \]
    However, $f(x) = (1 + x)^r$. Let $r = s_0 + s_1p^{f}+ \ldots + s_np^{nf}$ be the expansion of $r$ in base $p^{f}$. Using the identity $(1 + x)^{q} = (1 + x)$ in $R$, we see that $f(x) = (1 + x)^s$, where $s = s_0 + s_1 + \ldots +s_n$. If $r \geq q$, then $s < r$ and $s \equiv r \mod q - 1$. Iterating this process, we may assume that $s \in (0, q - 1]$. Since $s \equiv r \equiv 0 \mod q - 1$, we see that $s = q - 1$. Therefore $f(x) = (1 + x)^{q - 1}$ in $R$.
    
    If $i \in (0, q - 1)$, the coefficient of $x^i$ in $f(x)$ is
    \[
        {(p - 1)p^{f - 1} + (p - 1)p^{f - 2} + \cdots + (p - 1) \choose i_{f - 1}p^{f - 1} + i_{f - 2}p^{f - 2} + \cdots + i_0} = \prod_{j = 0}^{f - 1}{p - 1 \choose i_j} = (-1)^i
    \]
    in $\Fp$. If $i = 0$, however, the relation $x^{q - 1} - 1$ makes the leading term constant. So, the coefficient of $x^i$ receives contributions from the constant term and the leading term, which proves the lemma.
\end{proof}
We state the following theorem, which is the $m=0$ case in \cite[Theorem 2.18]{GJ25}.
\begin{theorem}\label{GJ Theorem 2.18}
Let $r = r_0 +  r_1p+ \ldots + r_{f-1}p^{f - 1}$ be such that $r_{j} \geq q=p^{f}$ for all $0 \leq j \leq f-1$, and $r$ is divisible by $q-1$. The map
        \[
            \psi : \otimes_{j=0}^{f-1} V_{r_{j}}^{\fr^{j}} \to \operatorname{ind}_{{B}(\Fq)}^{\GL_2(\Fq)}{\bf{1}} 
        \]
        that sends a polynomial $P=P_{0}(X_0,Y_0) \otimes \cdots \otimes P_{f-1}(X_{f-1},Y_{f-1}) \in \otimes_{j=0}^{f-1} V_{r_{j}}^{\fr^{j}} $ to the function $\psi_P : \GL_2(\Fq) \to \Fq$ defined by
        \[
            \psi_P\left(\begin{pmatrix}a & b \\ c & d\end{pmatrix}\right) = P_{0}(c,d) \otimes \cdots \otimes P_{f-1}(c^{p^{f-1}},d^{p^{f-1}})
        \]
        induces an isomorphism of $\GL_2(\Fq)$-representations
        \[
            \frac{\otimes_{j=0}^{f-1} V_{r_{j}}^{\fr^{j}} }{\langle\theta_0, \ldots, \theta_{f - 1}\rangle } \xrightarrow{\sim} \operatorname{ind}_{{B}(\Fq)}^{\GL_2(\Fq)} {\bf{1}}.
        \]
\end{theorem}
For brevity, let us define the notations. For $r = r_0 + \cdots + r_{f - 1}p^{f - 1}$, we write $\underline{X}^{\underline{r}}$ and $\underline{Y}^{\underline{r}}$ to represent $X_0^{r_0} \cdots X_{f - 1}^{r_{f - 1}}$ and $Y_0^{r_0} \cdots Y_{f - 1}^{r_{f - 1}}$, respectively. Also, 
$\underline{X}^{\underline{p-1}} \coloneq X_0^{p-1} \cdots X_{f - 1}^{p-1}$, $\underline{Y}^{\underline{p-1}} \coloneq Y_0^{p-1} \cdots Y_{f - 1}^{p-1}$ and $\underline{X}^{\underline{r-p +1}}\underline{Y}^{\underline{p - 1}} \coloneq  X_0^{r_0-(p-1)}Y_0^{p-1}\cdots X_{f - 1}^{r_{f - 1}-(p-1)} Y_{f - 1}^{p-1}$.
        
\begin{prop}\label{P1}
    Let $f \geq 2$. Let $r = r_0 + r_1 p + \ldots + r_{f - 1}p^{f - 1}$ be such that $r_j \geq q = p^f$ for all $0 \leq j \leq f - 1$. Assume $r$ is divisible by $q-1$. Then, the maps
    \begin{eqnarray*}
        \iota_0 : \otimes_{j=0}^{f-1} V_{0}^{\fr^{j}}  & \to & \frac{\otimes_{j=0}^{f-1} V_{r_{j}}^{\fr^{j}} }{\langle\theta_0, \ldots, \theta_{f - 1}\rangle } \text{ defined by}\\
        1 & \mapsto & \underline{X}^{\underline{r}} - \underline{X}^{\underline{r-p +1}}\underline{Y}^{\underline{p - 1}} + \underline{Y}^{\underline{r}}
    \end{eqnarray*}
    and 
    \begin{eqnarray*}
       \qquad \qquad \qquad \qquad \>\> \iota_{p - 1}: \otimes_{j=0}^{f-1} V_{p-1}^{\fr^{j}}  & \to & \frac{\otimes_{j=0}^{f-1} V_{r_{j}}^{\fr^{j}} }{\langle\theta_0, \ldots, \theta_{f - 1}\rangle } \text{ defined by} \\
        \underline{X}^{\underline{p-1-i}}\underline{Y}^{\underline{i}} & \mapsto & \underline{X}^{\underline{r-i}}\underline{Y}^{\underline{i}} \text{ if } (i_0, \ldots, i_{f - 1}) \neq (p - 1, \ldots, p - 1), \\
        \underline{Y}^{\underline{p-1}} & \mapsto & \underline{Y}^{\underline{r}} \text{ if } (i_0, \ldots, i_{f - 1}) = (p - 1, \ldots, p - 1)
    \end{eqnarray*}    
    induce an isomorphism
    \[
        \otimes_{j=0}^{f-1} V_{0}^{\fr^{j}}  \oplus \otimes_{j=0}^{f-1} V_{p-1}^{\fr^{j}}  \simeq \frac{\otimes_{j=0}^{f-1} V_{r_{j}}^{\fr^{j}} }{\langle\theta_0, \ldots, \theta_{f - 1}\rangle}.
    \]
\end{prop}
\begin{proof}
By Lemma \ref{GJ Theorem 2.18}, we have that $$\frac{\otimes_{j=0}^{f-1} V_{r_{j}}^{\fr^{j}} }{\langle\theta_0, \ldots, \theta_{f - 1}\rangle } \simeq \operatorname{ind}_{B(\Fq)}^{\GL_2(\Fq)} {\bf{1}} .$$ There is a copy of $\otimes_{j=0}^{f-1} V_{0}^{\fr^{j}} $ in $\operatorname{ind}_{B(\Fq)}^{\GL_2(\Fq)} {\bf{1}}$ given by the constant functions.
A brief calculation shows that under the isomorphism presented above, the polynomial $\underline{X}^{\underline{r}} - \underline{X}^{\underline{r-p +1}}\underline{Y}^{\underline{p - 1}} + \underline{Y}^{\underline{r}}$ maps to the constant function 1.
 
 Next, we prove that the map
    \begin{eqnarray*}
       \otimes_{j=0}^{f-1} V_{p-1}^{\fr^{j}}  & \to & \frac{\otimes_{j=0}^{f-1} V_{r_{j}}^{\fr^{j}} }{\langle\theta_0, \ldots, \theta_{f - 1}\rangle }\text{ defined by} \\
        \underline{X}^{\underline{p-1-i}}\underline{Y}^{\underline{i}} & \mapsto & \underline{X}^{\underline{r-i}}\underline{Y}^{\underline{i}} \text{ if } (i_0, \ldots, i_{f - 1}) \neq (p - 1, \ldots, p - 1), \\
        \underline{Y}^{\underline{p-1}} & \mapsto & \underline{Y}^{\underline{r}} \text{ if } (i_0, \ldots, i_{f - 1}) = (p - 1, \ldots, p - 1)
    \end{eqnarray*}
    is $\GL_2(\Fq)$-equivariant.

    Since $\otimes_{j=0}^{f-1} V_{p-1}^{\fr^{j}} $ is irreducible, it is generated by any non-zero vector. In particular, it is generated by $\underline{X}^{\underline{p-1}}$. Therefore to show that $\iota_{p - 1}$ is $\GL_2(\Fq)$-equivariant, it is enough to show that
    \[
        \iota_{p - 1}(g \underline{X}^{\underline{p-1}}) = g \iota_{p - 1}(\underline{X}^{\underline{p-1}}) \text{ for all } g \in \GL_2(\Fq).
    \]
    Noting the Bruhat decomposition $\GL_2(\Fq) = B \sqcup BwB$, it is enough to show
    \begin{eqnarray*}
        \iota_{p - 1}(M \underline{X}^{\underline{p-1}}) & = & M \iota_{p - 1}(\underline{X}^{\underline{p-1}}) \\
        \iota_{p - 1}(M \underline{Y}^{\underline{p-1}}) & = & M \iota_{p - 1}(\underline{Y}^{\underline{p-1}})
    \end{eqnarray*}
    for all $M \in B$. The first identity is easy to check since $B$ stabilizes $\underline{X}^{\underline{p-1}}$ and $\underline{X}^{\underline{r}}$. We now verify the second identity. Let $M = \begin{pmatrix}a & b \\ 0 & d\end{pmatrix}$. Then the left-hand side of the second equation above is
    \begin{eqnarray} \label{LHS of second eqn}
        && \iota_{p - 1}[(bX_0 + dY_0)^{p - 1} \cdots (b^{p^{f - 1}}X_{f - 1} + d^{p^{f - 1}}Y_{f - 1})^{p - 1}] \nonumber \\
        && = \iota_{p - 1}\left[\sum_{(i_0, \ldots, i_{f - 1}) = (0, \cdots, 0)}^{(p - 1, \ldots, p - 1)} (-1)^{i} b^{q - 1 - i}d^{i} X_0^{p - 1 - i_0}Y_0^{i_0} \cdots X_{f - 1}^{p - 1 - i_{f - 1}}Y_{f - 1}^{i_{f - 1}}\right] \nonumber \\
        && = \left[\sum_{\substack{(i_0, \ldots, i_{f - 1}) = (0, \cdots, 0) \\ (i_0, \ldots, i_{f - 1}) \neq (p - 1, \ldots, p - 1)}}^{(p - 1, \ldots, p - 1)} \!\!\!\!(-1)^{i} b^{q - 1 - i}d^{i} X_0^{r_0 - i_0}Y_0^{i_0} \cdots X_{f - 1}^{r_{f - 1} - i_{f - 1}}Y_{f - 1}^{i_{f - 1}}\right] + Y_{0}^{r_0} \cdots Y_{f - 1}^{r_{f - 1}}.
    \end{eqnarray}
    Next, we compute the right-hand side
    \begin{eqnarray}\label{RHS of second eqn}
        && M \iota_{p - 1} (\underline{Y}^{\underline{p-1}}) \nonumber \\
        && = M (\underline{Y}^{\underline{r}}) \nonumber \\
        && = (bX_0 + dY_0)^{r_0} \cdots (b^{p^{f - 1}}X_{f - 1} + d^{p^{f - 1}}Y_{f - 1})^{r_{f - 1}} \nonumber \\
        && = \sum_{(k_0, \ldots, k_{f - 1}) = (0, \ldots, 0)}^{(r_0, \ldots, r_{f - 1})} {r_0 \choose k_0} \cdots {r_{f - 1} \choose k_{f - 1}} b^{r - k} d^{k} X_{0}^{r_0 - k_0}Y_0^{k_0} \cdots X_{f - 1}^{r_{f - 1} - k_{f - 1}}Y_{f - 1}^{k_{f - 1}}.
    \end{eqnarray}
    
    Note that the coefficients of $\underline{X}^{\underline{r}}$ in \eqref{LHS of second eqn} and \eqref{RHS of second eqn} match, and so do the coefficients of $\underline{Y}^{\underline{r}}$. To show that the remaining expressions are equal, it is enough to check
    \[
        \sum_{\substack{(k_0, \ldots, k_{f - 1}) = (0, \ldots, 0) \\ k \equiv i \!\! \mod q - 1}}^{(r_0, \ldots, r_{f - 1})} {r_0 \choose k_0}\cdots {r_{f - 1} \choose k_{f - 1}} \equiv (-1)^i \mod p \text{ for }i \in (0,\  p^{f-1} - 1)
    \] and \[
        \sum_{\substack{(k_0, \ldots, k_{f - 1}) = (0, \ldots, 0) \\ (k_0, \ldots, k_{f - 1}) \neq (0, \ldots, 0) \\ (k_0, \ldots, k_{f - 1}) \neq (r_0, \ldots, r_{f - 1})  \\ k \equiv i \!\! \mod q - 1}  }^{(r_0, \ldots, r_{f - 1})} {r_0 \choose k_0}\cdots {r_{f - 1} \choose k_{f - 1}} \equiv 0 \mod p  \text{ for } i = 0.
    \]
    This is proved in Lemma~\ref{Sum of products of binomial coeffs is a sign}. 
    Finally, a dimension count shows $\otimes_{j=0}^{f-1} V_{0}^{\fr^{j}}  \oplus \otimes_{j=0}^{f-1} V_{p-1}^{\fr^{j}}  \simeq \frac{\otimes_{j=0}^{f-1} V_{r_{j}}^{\fr^{j}} }{\langle\theta_0, \ldots, \theta_{f - 1}\rangle }$.
\end{proof}
    
\begin{theorem}\label{main theorem}
    Let $f \geq 2$. Let $r = r_0 + r_1 p + \ldots + r_{f - 1}p^{f - 1}$ be such that $r_j \geq q = p^f$ for all $0 \leq j \leq f - 1$. Assume $r$ is divisible by $q - 1$. The map 
$$\operatorname{ind}_{IZ}^{G} {\bf{1}} \rightarrow \operatorname{ind}_{KZ}^{G} \frac{\otimes_{j=0}^{f-1} V_{r_{j}}^{\fr^{j}} }{\langle\theta_0, \ldots, \theta_{f - 1}\rangle } $$
$$\left[[\id, 1 \right]] \mapsto  \left[\id,  \underline{Y}^{\underline{r}} -\underline{X}^{\underline{r-p+1}}\underline{Y}^{\underline{p - 1}}\right]$$ induces isomorphisms
$$\frac{\operatorname{ind}_{IZ}^{G} {\bf{1}} }{\left(T_{1,2}T_{1,0}\right)} \simeq \operatorname{ind}_{KZ}^{G} \otimes_{j=0}^{f-1} V_{0}^{\fr^{j}}  \quad \text{and} \quad  \frac{\operatorname{ind}_{IZ}^{G} {\bf{1}} }{\left(\Id+T_{1,2}T_{1,0}\right)} \simeq \operatorname{ind}_{KZ}^{G} \otimes_{j=0}^{f-1} V_{p-1}^{\fr^{j}} .$$
Moreover, the operators $T_{-1,0} + T_{1,0}$ and $T_{-1,0}$ on the left correspond to the operator $T$ on the right under the isomorphism mentioned above.

\end{theorem}
\begin{proof}
    We first prove that the map 
    $$\operatorname{ind}_{IZ}^{G}{\bf{1}} \rightarrow \operatorname{ind}_{KZ}^{G} \frac{\otimes_{j=0}^{f-1} V_{r_{j}}^{\fr^{j}} }{\langle\theta_0, \ldots, \theta_{f - 1}\rangle } $$
$$\left[[\id, 1 \right]] \mapsto  \left[\id, \underline{Y}^{\underline{r}} -\underline{X}^{\underline{r-p + 1}}\underline{Y}^{\underline{p - 1}}\right]$$ is an isomorphism. Since $IZ \subseteq KZ \subseteq G$ with $ IZ \subseteq KZ$ of finite index, $\operatorname{ind}_{IZ}^{G} {\bf{1}}$ is isomorphic to $\operatorname{ind}_{KZ}^{G} \operatorname{ind}_{IZ}^{KZ}{\bf{1}}$ by sending 
$$\left[[\id, 1\right]] \mapsto \left[\id,\left[[[\id,1]]\right]\right],$$ 
A brief computation using Lemma \ref{GJ Theorem 2.18} shows that the map sending $\underline{Y}^{\underline{r}} -\underline{X}^{\underline{r-p +1}}\underline{Y}^{\underline{p - 1}}$ to $\left[[[\id,1]]\right]$ is an isomorphism between $ \frac{\otimes_{j=0}^{f-1} V_{r_{j}}^{\fr^{j}}}{\langle\theta_0, \ldots, \theta_{f - 1}\rangle }$ and $\operatorname{ind}_{IZ}^{KZ} {\bf{1}}$.\\
 Since $-T_{1,2}T_{1,0}$ is an idempotent in the non-commutative Hecke algebra, we get a splitting 
 $$\operatorname{ind}_{IZ}^{G} {\bf{1}}={\im}(-T_{1,2}T_{1,0}) \oplus {\im}(\Id+T_{1,2}T_{1,0}).$$
 We also know from Proposition \ref{P1} that 
 $$\operatorname{ind}_{KZ}^{G} \frac{\otimes_{j=0}^{f-1} V_{r_{j}}^{\fr^{j}} }{\langle\theta_0, \ldots, \theta_{f - 1}\rangle } = \operatorname{ind}_{KZ}^{G} \otimes_{j=0}^{f-1} V_{0}^{\fr^{j}}  \oplus \operatorname{ind}_{KZ}^{G} \otimes_{j=0}^{f-1} V_{p-1}^{\fr^{j}} $$
 where $\operatorname{ind}_{KZ}^{G} \left(\otimes_{j=0}^{f-1} V_{p-1}^{\fr^{j}}\right)$ is identified with the subspace generated by $\left[\id, \underline{X}^{\underline{r}} \right]$ and $\operatorname{ind}_{KZ}^{G}\left(\otimes_{j=0}^{f-1} V_{0}^{\fr^{j}}\right)$ is identified with subspace generated by $\left[\id, \underline{X}^{\underline{r}} - \underline{X}^{\underline{r-p+ 1}}\underline{Y}^{\underline{p - 1}} + \underline{Y}^{\underline{r}} \right]$. 
 Therefore, to prove the first part of this theorem, we need to show that 
  \begin{itemize}
      \item ${\im}(-T_{1,2}T_{1,0}) \mapsto \operatorname{ind}_{KZ}^{G} (\otimes_{j=0}^{f-1} V_{p-1}^{\fr^{j}} )$,
      \item ${\im}(\Id+T_{1,2}T_{1,0}) \mapsto \operatorname{ind}_{KZ}^{G} (\otimes_{j=0}^{f-1} V_{0}^{\fr^{j}} )$.
      \end{itemize} 
      Indeed, under the map $\operatorname{ind}_{IZ}^{G} {\bf{1}} \rightarrow \operatorname{ind}_{KZ}^{G} \frac{\otimes_{j=0}^{f-1} V_{r_{j}}^{\fr^{j}} }{\langle\theta_0, \ldots, \theta_{f - 1}\rangle}$, the element $T_{1,2}\left[[\id,1\right]]$ maps to 
      \[
                \sum_{\lambda \in I_1}\left[\begin{pmatrix}1 & 0 \\ \varpi \lambda & \varpi\end{pmatrix}, \underline{Y}^{\underline{r}} -\underline{X}^{\underline{r-p+1}}\underline{Y}^{\underline{p - 1}}\right] = \sum_{\lambda \in I_1}\left[\beta\begin{pmatrix}1 & \lambda \\ 0 & 1\end{pmatrix}w, \underline{Y}^{\underline{r}} -\underline{X}^{\underline{r-p +1}}\underline{Y}^{\underline{p - 1}}\right].
            \]
            Since $\begin{pmatrix}1 & \lambda \\ 0 & 1\end{pmatrix}w \in KZ$, we can move it to the other side to get 
            \[
                \sum_{\lambda \in I_1}\left[\beta, X_0^{r_0} \cdots X_{f - 1}^{r_{f - 1}}- (\lambda X_0+Y_0)^{r_0 -(p - 1)}X_0^{p - 1} \cdots (\lambda^{p^{f - 1}} X_{f - 1}+Y_{f-1})^{r_{f - 1} - (p - 1)}X_{f - 1}^{p - 1}\right].
            \]
            Expanding $(\lambda^{p^{j}} X_j + Y_j)^{r_j-(p - 1)}$ for all $0 \leq j \leq f-1$, above expression is 
            
\begin{eqnarray*}
       =- \left[\beta,  \sum_{(k_0, \ldots, k_{f - 1}) = (0, \ldots, 0)}^{(r_0-(p-1), \ldots, r_{f - 1}-(p-1))} {r_0-(p-1) \choose k_0} \cdots {r_{f - 1}-(p-1) \choose k_{f - 1}} \sum_{\lambda \in \Fq^{\times}}\lambda^{r-(q-1)-k}  \underline{Y}^{\underline{r-k}}\underline{X}^{\underline{k}}\right]
            \end{eqnarray*}
            $\hspace{0.4 cm}-\left[\beta, \underline{Y}^{\underline{r-p+1}}\underline{X}^{\underline{p-1}}\right]$.\\
        Using certain congruences between monomials in $V_r$ modulo $\langle \theta_0, \ldots, \theta_{f - 1}\rangle$ and summing over $\lambda$ using the identity 
            \begin{eqnarray}\label{sum of powers of roots of unity}
            \sum_{\mu \in \Fq}\mu^l = 
            \begin{cases}
                0 & \text{ if } q - 1 \nmid l \\
                q - 1 & \text{ if } q - 1 \mid l,
                \end{cases}
                & \text{ for $l \geq 1$,}
            \end{eqnarray} 
             and Lemma~\ref{Sum of products of binomial coeffs is a sign} with $r_j$ replaced with $r_j - (p - 1)$,
            we see that
            \begin{eqnarray}\label{Image of T12 maps to}
                T_{1, 2}[[\id, 1]] \mapsto \left[\beta,\underline{X}^{\underline{r}}\right].
            \end{eqnarray}
\begin{itemize}
\item Using \eqref{Image of T12 maps to}, we see that
            $$T_{1, 2}T_{1, 0}[[\id, 1]] = T_{1, 2}[[\beta, 1]] \mapsto \left[\id, \underline{X}^{\underline{r}}\right].
            $$ 
            Moreover, $\left[\id, \underline{X}^{\underline{r}}\right]$ generates $\operatorname{ind}_{KZ}^{G} \left(\otimes_{j=0}^{f-1} V_{p-1}^{\fr^{j}} \right)$ in $\operatorname{ind}_{KZ}^{G} \frac{\otimes_{j=0}^{f-1} V_{r_{j}}^{\fr^{j}} }{\langle\theta_0, \ldots, \theta_{f - 1}\rangle }$.
\item Again using \eqref{Image of T12 maps to}, we see that $$(\Id+T_{1, 2}T_{1, 0})[[\id, 1]] = [[\id, 1]] + T_{1, 2}[[\beta, 1]]$$
$$[[\id, 1]] + T_{1, 2}[[\beta, 1]] \mapsto \left[\id,  \underline{X}^{\underline{r}} - \underline{X}^{\underline{r-p +1}}\underline{Y}^{\underline{p - 1}} + \underline{Y}^{\underline{r}}\right]$$
Moreover $\left[\id,  \underline{X}^{\underline{r}} - \underline{X}^{\underline{r-p+1}}\underline{Y}^{\underline{p - 1}} + \underline{Y}^{\underline{r}}\right]$ generates $\operatorname{ind}_{KZ}^{G} \left(\otimes_{j=0}^{f-1} V_{0}^{\fr^{j}} \right)$.
        \end{itemize}
        We have proved that the natural map induces isomorphisms
        \begin{eqnarray}\label{Main comparison theorems in the non-commutative case first map}
            \frac{\operatorname{ind}_{IZ}^{G}{\bf{1}}}{(T_{1, 2}T_{1, 0})} \xrightarrow{\sim} \operatorname{ind}_{KZ}^{G} \!\left(\otimes_{j=0}^{f-1} V_{0}^{\fr^{j}} \right)\!\! \text{  and } \frac{\operatorname{ind}_{IZ}^{G}{\bf{1}}}{(\Id + T_{1, 2}T_{1, 0})} \xrightarrow{\sim} \operatorname{ind}_{KZ}^{G} \!\left(\otimes_{j=0}^{f-1} V_{p-1}^{\fr^{j}} \right).\!\!
        \end{eqnarray}
            These maps are respectively given by
            \[
                [[\id, 1]] \mapsto [\id, 1] \quad \text{ and } \quad [[\id, 1]] \mapsto \left[\id,-(\underline{X}^{\underline{p - 1}})\right].
            \]
            Indeed, the first map sends $[[\id, 1]]$ to \begin{eqnarray*}
                \left[\id, \underline{Y}^{\underline{r}} -\underline{X}^{\underline{r-p +1}}\underline{Y}^{\underline{p - 1}}\right] \equiv \left[\id, \underline{X}^{\underline{r}} - \underline{X}^{\underline{r-p +1}}\underline{Y}^{\underline{p - 1}} + \underline{Y}^{\underline{r}}\right] \mod \operatorname{ind}_{KZ}^{G} \left(\otimes_{j=0}^{f-1} V_{p-1}^{\fr^{j}} \right)
            \end{eqnarray*} and the second map sends $[[\id, 1]]$ to 
            \begin{eqnarray*}
                \left[\id, \underline{Y}^{\underline{r}} -\underline{X}^{\underline{r-p +1}}\underline{Y}^{\underline{p - 1}}\right] \equiv \left[\id,-(\underline{X}^{\underline{r}})\right] \mod \operatorname{ind}_{KZ}^{G} \left(\otimes_{j=0}^{f-1} V_{0}^{\fr^{j}} \right).
            \end{eqnarray*}
            Now using \eqref{relations in non-reqular case} and \eqref{T-10 in the non-commutative case}, we see that ${\im} T_{1, 2}T_{1, 0} \subseteq {\Ker} (T_{-1,0} + T_{1,0})$ and ${\im} (1 + T_{1, 2}T_{1, 0}) \subseteq {\Ker} T_{-1, 0}$. Therefore, they are operators on the quotients
            \[
                \frac{\operatorname{ind}_{IZ}^{G}{\bf{1}}}{(T_{1, 2}T_{1, 0})} \quad \text{ and } \quad \frac{\operatorname{ind}_{IZ}^{G}
                {\bf{1}}}{(\Id + T_{1, 2}T_{1, 0})},
            \]
            respectively. Then, under the first map in \eqref{Main comparison theorems in the non-commutative case first map}, we see that
            \[
                (T_{-1, 0} + T_{1, 0})[[\id, 1]] = \sum_{\lambda \in I_1}\left[\left[\begin{pmatrix}\varpi & \lambda \\ 0 & 1\end{pmatrix}, 1\right]\right] + [[\beta, 1]] \mapsto \sum_{\lambda \in I_1}\left[\begin{pmatrix}\varpi & \lambda \\ 0 & 1\end{pmatrix}, 1\right] + \left[\alpha, 1\right] = T[\id, 1]
            \]
            since $\beta = \alpha w$ and $w \in KZ$. Similarly, under the second map in \eqref{Main comparison theorems in the non-commutative case first map}, we get
            \[
                T_{-1, 0}[[\id, 1]] = \sum_{\lambda \in I_1}\left[\left[\begin{pmatrix}\varpi & \lambda \\ 0 & 1\end{pmatrix}, 1\right]\right] \mapsto \sum_{\lambda \in I_1}\left[\begin{pmatrix}\varpi & \lambda \\ 0 & 1\end{pmatrix}, -(\underline{X}^{\underline{p - 1}})\right] = T[\id, -(\underline{X}^{\underline{p-1}}) ]. \qedhere
            \]
\end{proof}

\section{Self Extensions of a Supersingular representation of $\rm{GL_2}(\mathbb{Q}_p)$}
Throughout this section, let $p$ be an odd prime. Let $G=\rm{GL_2}(\mathbb{Q}_p)$ and $K=\rm{GL_2}(\mathbb{Z}_p)$ be its standard maximal compact subgroup. Let $I$ be the Iwahori subgroup of $G$. Let $\pi_{r}=\frac{\operatorname{ind}_{KZ}^{G} V_r}{(T)}$ where, $r\in \{0,p-1\}$ be a supersingular representation of $G$.

Breuil \cite{Bre03} has given the following isomorphism
$$\frac{\operatorname{ind}_{KZ}^{G} V_0}{(T)} \cong \frac{\operatorname{ind}_{KZ}^{G} V_{p-1}}{(T)}.$$
Before we come to results, we need to consider the implications of the above isomorphism for finding the self-extensions of $\pi_r$. A key property of the Hecke operator $T$ is its injectivity
which implies that we have an exact sequence $${\displaystyle 0\longrightarrow \pi_r{\longrightarrow{}\frac{\operatorname{ind}_{KZ}^{G}V_r}{(T^2)}}{\longrightarrow{}\pi_r \longrightarrow 0}}.$$

The following lemma is well-known. We could not find its proof in the literature. Hence, we provide a short proof. One can compare this lemma on the automorphic side with \cite[Theorem 7.53]{BGR18} on the Galois side using the non-semisimple mod $p$ local Langlands correspondence.

\begin{lemma}
   For $r\in \{0,p-1\}$, the exact sequence 
    $${\displaystyle 0\longrightarrow \pi_r{\overset{T}{\longrightarrow} \frac{\operatorname{ind}_{KZ}^{G}V_{0}}{(T^2)}}{\longrightarrow{}\pi_r \longrightarrow 0}}$$
    is non-split.
\end{lemma}
\begin{proof}
Since $\pi_0 \cong \pi_{p-1}$, it suffices to prove that the sequence is non-split for $r=0$. Suppose that the given sequence is split. So, there exists a retraction map $S : \frac{\KZind V_0}{(T^2)} \rightarrow \pi_0$ for the map $T$.
Using the fact that $[\id, 1]$ is fixed under $I(1)$, we see that its image in $\frac{\KZind V_0}{(T^2)}$ is of the form $[\id, a] + [\alpha, b]$ in $\pi_0$ for some $a, b \in \Fp$. Also, since $[\id, 1]$ is fixed by $w$ in $\frac{\KZind V_0}{(T^2)}$, its image is fixed by $w$. Then, a quick computation using \cite[Lemma 20]{BL94} shows $b = 0$.
Writing $T[\id, 1]$ as a sum of $\GL_2(\Q_p)$-translates of $[\id, 1]$, we see that $S$ commutes with $T$. Consequently, $S(T[\id, 1])=0$ in $\pi_0$. However, since $S$ is a retraction map for $T$, we get a contradiction.
\end{proof}

It is known in \cite{Pas10} and \cite{BP12} that the one dimensional subspace of the three dimensional $\Ext_{G}^{1}\left(\pi_r, \pi_r \right)$ spanned by the class of this self-extension consists of all the self extensions of $\pi_r$ which have a four-dimensional space of $I(1)$-invariants, where $I(1)$ denotes the pro-$p$ Iwahori subgroup of $G$, and $\Ext$ is computed in the category of smooth representations with a given central character.

In this section, we focus on giving a realization of the self-extensions of $\pi_r$, belonging to the above-mentioned one-dimensional space, using the Iwahori-Hecke model for the \emph{non-regular} case, i.e. $r \in \{0,p-1\}$, given in the recent work of the first author \cite{Chi25}. For the \emph{regular} case $0<r<p-1$, see \cite{AB15}.

The following theorem gives the Iwahori-Hecke model for a supersingular representation of $\rm{GL_2}(\mathbb{Q}_p)$ in the \emph{non-regular} case.

\begin{theorem}
 We have, $$\frac{\operatorname{ind}_{IZ}^{G} \bf{1}}{(T_{1,2}T_{1,0})} \cong \operatorname{ind}_{KZ}^{G} V_0  \  \text{ and  } \ \frac{\operatorname{ind}_{IZ}^{G} \bf{1}}{(\Id+T_{1,2}T_{1,0})} \cong \operatorname{ind}_{KZ}^{G}V_{p-1} .$$ Moreover, the operators $T_{-1,0}+T_{1,0}$ and $T_{-1,0}$ on the left correspond to the operator $T$ on the right.
\end{theorem}
\begin{proof}
See Theorem 3.13 in \cite{Chi25}. 
\end{proof}

Explicitly, Breuil's isomorphism in the \emph{non-regular} case is given as follows.
\begin{theorem}\textbf{(Breuil's isomorphism)}
   We have, 
   $$\frac{\operatorname{ind}_{IZ}^{G} \bf{1}}{(T_{-1,0})+(T_{1,2}+T_{1,0})} \cong \frac{\operatorname{ind}_{IZ}^{G} \bf{1}}{(T_{1,0}+T_{-1,0})+(T_{1,2})}.$$ 
\end{theorem}
\begin{proof}
Indeed, $T_{1,0}$ induces an isomorphism.   
\end{proof}

\begin{theorem}\label{T3}
       Let $\pi_{0}=\frac{\operatorname{ind}_{KZ}^{G} V_0}{(T)}$ be an irreducible supersingular representation of $\rm{GL_2}(\mathbb{Q}_p)$. Then, a basis for the subspace of $\Ext_{G}^{1}\left(\pi_{0}, \pi_{0} \right)$ consisting of isomorphism class of short exact sequence 
    $${\displaystyle 0\longrightarrow \pi_{0}{\longrightarrow{}\tau}{\longrightarrow{}\pi_{0} \longrightarrow 0}},$$ where the representation $\tau$ has a four dimensional space of $I(1)$-invariants, is given by isomorphism class of 
    $${\displaystyle 0\longrightarrow \pi_{0}{\longrightarrow{}\frac{\operatorname{ind}_{IZ}^{G} \bf{1}}{(T_{-1,0})+(T_{1,2}+T_{1,0})^2} \cong \frac{\operatorname{ind}_{KZ}^{G} V_0}{(T^2)}}{\longrightarrow{}\pi_{0} \longrightarrow 0}}.$$
\end{theorem}

\begin{theorem}\label{T4}
    Let $\pi_{p-1}=\frac{\operatorname{ind}_{KZ}^{G}V_{p-1}}{(T)}$ be an irreducible supersingular representation of $\rm{GL_2}(\mathbb{Q}_p)$. Then, a basis for the subspace of $\Ext_{G}^{1}\left(\pi_{p-1}, \pi_{p-1} \right)$ consisting of isomorphism class of short exact sequence 
    $${\displaystyle 0\longrightarrow \pi_{p-1}{\longrightarrow{}\tau'}{\longrightarrow{}\pi_{p-1} \longrightarrow 0}},$$ where the representation $\tau'$ has a four dimensional space of $I(1)$-invariants, is given by isomorphism class of 
    $${\displaystyle 0\longrightarrow \pi_{p-1}{\longrightarrow{}\frac{\operatorname{ind}_{IZ}^{G} \bf{1}}{(T_{1,0}+T_{-1,0})+(T^2_{1,2})} \cong \frac{\operatorname{ind}_{KZ}^{G} V_{p-1}}{(T^2)}}{\longrightarrow{}\pi_{p-1} \longrightarrow 0}}.$$

\end{theorem}

Before we prove the theorems we note that $\tau=\frac{\operatorname{ind}_{IZ}^{G} \bf{1}}{(T_{-1,0})+(T_{1,2}+T_{1,0})^2}$ and $\tau'=\frac{\operatorname{ind}_{IZ}^{G} \bf{1}}{(T_{1,0}+T_{-1,0})+(T^2_{1,2})}$ are isomorphic by 
defining a map 
$$\Phi : \frac{\operatorname{ind}_{IZ}^{G} \bf{1}}{X=[(T_{1,0}+T_{-1,0})+(T^2_{1,2})]} \longrightarrow \frac{\operatorname{ind}_{IZ}^{G} \bf{1}}{[(T_{-1,0})+(T_{1,2}+T_{1,0})^2]=Y} $$ as 
$$\Phi(f+X)=(-T_{1,2})(f)+Y.$$
It is easy to check that $\Phi$ is a well-defined $ G$-linear map. Also, $\Phi$ makes the diagram commute. Hence, $\Phi$ is an isomorphism.
\subsection{}{\textbf{Proof of Theorem \ref{T3}}}
\begin{proof}
We prove that $\tau=\frac{\operatorname{ind}_{IZ}^{G} \bf{1}}{(T_{-1,0})+(T_{1,2}+T_{1,0})^2} \in \Ext_{G}^{1}\left(\pi_{0}, \pi_{0} \right)$. Let $\pi_{0} \cong \frac{\operatorname{ind}_{IZ}^{G} \bf{1}}{(T_{-1,0})+(T_{1,2}+T_{1,0})}$ be an irreducible supersingular representation of $\rm{GL_2}(\mathbb{Q}_p)$.  Consider the natural surjection
$$\Psi: \tau \longrightarrow \frac{\operatorname{ind}_{IZ}^{G} \bf{1}}{(T_{-1,0})+(T_{1,2}+T_{1,0})} $$ 
$$f+(T_{-1,0})+(T_{1,2}+T_{1,0})^2 \mapsto f+(T_{-1,0})+(T_{1,2}+T_{1,0}).$$
We can easily see that ${\Ker}\Psi=\frac{(T_{-1,0})+(T_{1,2}+T_{1,0})}{Y=[(T_{-1,0})+(T_{1,2}+T_{1,0})^2]}$. We claim that ${\Ker}\Psi \cong \frac{\operatorname{ind}_{IZ}^{G} \bf{1}}{(T_{-1,0})+(T_{1,2}+T_{1,0})}$. Consider a map 
$$\phi: \operatorname{ind}_{IZ}^{G} {\bf{1}}  \longrightarrow \frac{(T_{-1,0})+(T_{1,2}+T_{1,0})}{Y=[(T_{-1,0})+(T_{1,2}+T_{1,0})^2]} $$
$$f \mapsto (T_{1,2}+T_{1,0})(f)+Y$$
Note that $\phi$ is a $ G$-linear onto map. We show that ${\Ker}\phi= (T_{-1,0})+(T_{1,2}+T_{1,0}) $. For that, we prove the following result.
\begin{lemma}\label{L2}
    We have, ${\Ker}(T_{-1,0})={\im}(T_{1,2}+T_{1,0})$ and ${\Ker}(T_{1,2}+T_{1,0})={\im}(T_{-1,0})$.
\end{lemma}
\begin{proof}
Since $ T_{1,0}$ is invertible, ${\Ker}(T_{1,2}T_{1,0})={\Ker}(T_{1,0}T_{1,2}T_{1,0})={\Ker}(T_{-1,0})$ and ${\im}(T_{1,2}+T_{1,0})={\im}(T_{1,2}T_{1,0}+\Id)$. Using  \cite[Proposition 3.10]{Chi25}, we get the first equality. For the second equality, a similar proof applies using ${\im}(T_{1,0}T_{1,2})={\Ker}(T_{1,0}T_{1,2}+\Id)$.  
\end{proof}

\begin{remark}
     The Lemma \ref{L2} is true for any local field $F$.
 \end{remark}

 We can observe from the above Lemma \ref{L2} that ${\im} (T_{1,2}+T_{1,0})+{\im} (T_{-1,0})\subseteq {\Ker}\phi$. Since $\frac{\operatorname{ind}_{IZ}^{G} \bf{1}}{(T_{-1,0})+(T_{1,2}+T_{1,0})}$ is irreducible and $\phi$ is nonzero, we have  ${\Ker}\phi= (T_{-1,0})+(T_{1,2}+T_{1,0})$. Hence,
$$\frac{\operatorname{ind}_{IZ}^{G} \bf{1}}{ (T_{-1,0})+(T_{1,2}+T_{1,0})} \cong \frac{(T_{-1,0})+(T_{1,2}+T_{1,0})}{X=[(T_{-1,0})+(T_{1,2}+T_{1,0})^2]}={\Ker}\Psi,$$
which proves our claim $\tau=\frac{\operatorname{ind}_{IZ}^{G} \bf{1}}{( T_{-1,0})+(T_{1,2}+T_{1,0})^2} \in \Ext_{G}^{1}\left(\pi_{0}, \pi_{0} \right)$.
 
\end{proof}

\subsection{}{\textbf{Proof of Theorem \ref{T4}}}

\begin{proof}
We prove that $\tau'=\frac{\operatorname{ind}_{IZ}^{G} \bf{1}}{(T_{1,0}+T_{-1,0})+(T^2_{1,2})} \in \Ext_{G}^{1}\left(\pi_{p-1}, \pi_{p-1} \right)$. Let $\pi_{p-1} \cong \frac{\operatorname{ind}_{IZ}^{G} \bf{1}}{(T_{1,0}+T_{-1,0})+(T_{1,2})}$ be an irreducible supersingular representation of $\rm{GL_2}(\mathbb{Q}_p)$.  Consider the natural surjection
$$\Psi': \tau' \longrightarrow \frac{\operatorname{ind}_{IZ}^{G} \bf{1}}{(T_{1,0}+T_{-1,0})+(T_{1,2})} $$ 
$$f+(T_{1,0}+T_{-1,0})+(T^2_{1,2}) \mapsto f+(T_{1,0}+T_{-1,0})+(T_{1,2}).$$
We can easily see that ${\Ker}\Psi'=\frac{(T_{1,0}+T_{-1,0})+(T_{1,2})}{X=[(T_{1,0}+T_{-1,0})+(T^2_{1,2})]}$. We claim that ${\Ker}\Psi' \cong \frac{\operatorname{ind}_{IZ}^{G} \bf{1}}{(T_{1,0}+T_{-1,0})+(T_{1,2})}$. Consider a map 
$$\phi': \operatorname{ind}_{IZ}^{G} {\bf{1}}  \longrightarrow \frac{(T_{1,0}+T_{-1,0})+(T_{1,2})}{X=[(T_{1,0}+T_{-1,0})+(T^2_{1,2})]} $$
$$f \mapsto (T_{1,2})(f)+X$$
Note that $\phi'$ is a $ G$-linear onto map. We show that ${\Ker}\phi'= (T_{1,0}+T_{-1,0})+(T_{1,2})$. For that, we prove the following result.
\begin{lemma}\label{L1}
    We have, ${\Ker}(T_{1,0}+T_{-1,0})={\im}(T_{1,2})$ and ${\Ker}(T_{1,2})={\im}(T_{1,0}+T_{-1,0})$.
\end{lemma}
\begin{proof}
Since $T_{1,0}$ is invertible, ${\Ker}(T_{1,0}+T_{-1,0})={\Ker}(\Id+T_{1,2}T_{1,0})={\im}(T_{1,2}T_{1,0})={\im}(T_{1,2})$. For the second equality, use ${\Ker}(T_{1,2}T_{1,0})={\im}(T_{1,2}T_{1,0}+\Id)$ and $T_{1,0}$ induces isomorphism  ${\Ker}(T_{1,2}) \cong {\Ker}(T_{1,2}T_{1,0})$ and ${\im}(T_{1,0}+T_{-1,0}) \cong {\im}(T_{1,2}T_{1,0}+\Id)$. 
\end{proof}

 \begin{remark}
     The Lemma \ref{L1} is true for any local field $F$.
 \end{remark}

We can observe from the above Lemma \ref{L1} that ${\im}(T_{-1,0}+T_{1,0})+{\im}(T_{1,2}) \subseteq {\Ker}\phi'$. Since $\frac{\operatorname{ind}_{IZ}^{G} \bf{1}}{(T_{1,0}+T_{-1,0})+(T_{1,2})}$ is irreducible and $\phi'$ is nonzero, we have ${\Ker}\phi'= (T_{1,0}+T_{-1,0})+(T_{1,2})$. Hence,
$$\frac{\operatorname{ind}_{IZ}^{G} \bf{1}}{(T_{1,0}+T_{-1,0})+(T_{1,2})} \cong \frac{(T_{1,0}+T_{-1,0})+(T_{1,2})}{X=[(T_{1,0}+T_{-1,0})+(T^2_{1,2})]}={\Ker}\Psi',$$
which proves our claim $\tau'=\frac{\operatorname{ind}_{IZ}^{G} \bf{1}}{(T_{1,0}+T_{-1,0})+(T^2_{1,2})} \in \Ext_{G}^{1}\left(\pi_{p-1}, \pi_{p-1} \right)$.
\end{proof}

\begin{remark}
 Since we have the following isomorphism (as vector spaces) 
  $$\left(\frac{\operatorname{ind}_{IZ}^{G} \bf{1}}{(T_{-1,0})+(T_{1,2}+T_{1,0})}\right)^{I(1)} \cong \left(\frac{\operatorname{ind}_{IZ}^{G} \bf{1}}{(T_{1,0}+T_{-1,0})+(T_{1,2})}\right)^{I(1)}=\left<\overline{[[\id,1]]}, \overline{[[\beta,1]]}\right>,$$
 
  \begin{equation*}
  \left(\frac{\operatorname{ind}_{IZ}^{G} \bf{1}}{(T_{-1,0})+(T_{1,2}+T_{1,0})^2}\right)^{I(1)}
  \end{equation*}
  \begin{equation*}
   = \left <\overline{[[\id,1]]}, \overline{[[\beta,1]]},\overline{(T_{1,2}+T_{1,0})([[\id,1]])},\overline{(T_{1,2}+T_{1,0})([[\beta,1]])}\right>, 
  \end{equation*}
  and 
  $$\left(\frac{\operatorname{ind}_{IZ}^{G} \bf{1}}{(T_{1,0}+T_{-1,0})+(T^2_{1,2})}\right)^{I(1)}=\left <\overline{[[\id,1]]}, \overline{[[\beta,1]]},\overline{T_{1,2}([[\id,1]])},\overline{T_{1,2}([[\beta,1]])}\right>,$$
    
\end{remark}

\begin{lemma}
We have, $(T_{1,2}+T_{1,0})^2=T_{1,0}\ (T_{-1,0}+T_{1,0})^2 \ T_{1,0}$.
\end{lemma}
\begin{proof}
Indeed, $T_{1,0}\ (T_{-1,0}+T_{1,0}) \ T_{1,0} \ T_{1,0}\ (T_{-1,0}+T_{1,0}) \ T_{1,0}=(T_{1,2}+T_{1,0})^2.$
\end{proof}

\vspace{1cm}
\noindent{\bf  Acknowledgements.}
We thank Prof. Eknath Ghate, Prof. U. K. Anandavardhanan, Prof. Gautam Borisagar, and Arindam Jana for fruitful discussions.
The first author gratefully acknowledges support from the HRI postdoctoral fellowship.
The second author thanks Prof. Eknath Ghate for the invitation and warm hospitality during his visit to TIFR, Mumbai. Part of the work was done during the second author's visit to TIFR. Additionally, the second author acknowledges the financial assistance as a Junior Research Fellowship (NET) from UGC, India.

\printbibliography

\end{document}